\documentclass[amsmath,article,aps,fleqn]{revtex4}

\begin{document}

\title{Finite, closed-form expressions for the partition function and for  Euler, Bernoulli, and Stirling numbers}
\author{Jerome Malenfant}
\affiliation{American Physical Society\\ Ridge, NY}
\date{\today}
\begin{abstract}  We find general solutions to the generating-function equation $\sum c_q^{(X)} z^q = F(z)^X$,
where $X$ is a complex number and $F(z)$ is a convergent power series with $F(0) \neq 0$.  
We then use these results to derive finite expressions containing only integers or simple fractions for partition 
functions and for Euler, Bernoulli, and Stirling numbers.  %SourceDoc 
\end{abstract}

\maketitle
\newtheorem{I}{Theorem}
\newtheorem{II}{Proposition}
\newtheorem {4}{Corollary}
\newtheorem {1}{Lemma}
\newtheorem {2}[1]{Lemma}
\newtheorem {3}[1]{Lemma}

\section{Introduction and main result}

Generating functions are often used as a compact way to define special number sequences and functions as the coefficients in a power-series expansion of more elementary functions.  One example of this is the partition function $p(n)$, (sequence A000041 in OEIS \cite {integer}), the number of partitions of $n$ into positive integers.  It has the generating function \cite{Partition}
	\begin{eqnarray}
         \sum_{q=0}^{\infty} p(q) z^q = \prod_{k=1}^{\infty} \frac{1}{1-z^k} = 1 +z + 2z^2 +3z^3 +5z^4 + \cdots
	\end{eqnarray}
Some other examples of  generating functions are those for the Bernoulli and the Euler numbers:
         \begin{subequations}
         \begin{eqnarray}
         \sum_{q=0}^{\infty} \frac{ B_q}{q!} ~z^q &=& \frac{z}{e^z-1} ,\\
         \sum _{q=0}^{\infty} \frac{E_{2q}}{(2q)!}~ z^{2q}&=&  \frac{1}{\cosh z}       
         \end{eqnarray}
         \end{subequations}
The solutions of generating-function equations can often be found from recursion relations.   Some of them can also be 
calculated directly by using various expressions, such as Laplace's determinental formula for the Bernoulli numbers \cite{Korn}, which we write in the form:
	\begin{eqnarray}
	B_n	=    n!  \left| \begin{array}{ccccc}
          1 & 0 & \cdots  &~~ 0 & ~~1 \\
          \frac{1}{2!}  & 1 & ~ & ~~ 0 & ~~ 0 \\
          \vdots  &  ~ & \ddots  & ~  & ~~\vdots \\  
          \frac{1}{n!} & \frac{1}{(n-1)!} &~ & ~~1  & ~~ 0 \\            
          \frac{1}{ (n+1)!}&\frac{1}{n!} & \cdots  & ~~ \frac{1}{2!}  &  ~~0
          \end{array} \right|
          \end{eqnarray}
Vella \cite{Vella} has derived expressions for $B_n$ and for $E_n$ as sums over the partitions (or, alternatively, over the 
 compositions) of $n$ using a method based on the Faa di Bruno formula for the higher derivatives of composite functions.  
 Concerning the partition function, Rademacher  \cite {Partition, Rade} derived an exact formula, an improvement over the 
 Hardy-Ramanujan asymptotic formula, which however involves an infinite sum of rather complicated, non-integer terms.   
More recently,  Bruinier and Ono \cite{Ono} have derived an explicit formula for the partition function as  a finite sum of algebraic numbers that requires finding a sufficiently precise approximation to an auxiliary function.
               
In the following theorem, we present methods for solving generating-function equations for the case where the right-hand side
is expressible in terms of a convergent, non-zero-near-the-origin power series.  We will then apply this theorem to solve for the partition function and other objects which have the appropriate generating functions.  In the following, the notations
	 \begin{eqnarray*}
	 \left( \begin{array}{c} n \\ q_1, \ldots , q_K \end{array} \right)
          \equiv \frac{ n!}{q_1! \cdots q_K!} ~ \delta_{n,q_1+ \cdots + q_K}
	 \end{eqnarray*}
and
          \begin{eqnarray*}          
	 \left( \begin{array} {c} X \\ N \end{array} \right)  \equiv \frac{X(X-1) \cdots (X-N+1)}{N!} 
	 \end{eqnarray*}
will denote, respectively,   multinomial coefficients and generalized binomial coefficients.  

Our main result is:
          \begin{I} Let $F(z)$ be a holomorphic function in a neighborhood of the origin with the series expansion
           $ \sum _{q=0}^{\infty}    a_qz^q$,  with $a_0 \neq 0$, and let $X$ be a complex number.
          Then the coefficients $c_p^{(X)}$ in the generating-function equation
          \begin{eqnarray*}         
          \sum_{p=0}^{\infty} c^{(X)}_pz^p  = F(z)^X     
          \end{eqnarray*}    
are given by the two equivalent expressions:
          \begin{eqnarray*}
          ({\rm I}) ~~~~          c_p^{(X)} &=&    \sum_{ 0 \leq k_1, \ldots, k_p \leq p} 
          \left( \begin{array} {c} X \\ K \end{array} \right) 
          \left( \begin{array}{c} K \\ k_1, \ldots, k_p \end{array} \right) 
          \delta_{p, \sum mk_m} ~a_0^{X-K}  a_1^{k_1} \cdots a_p^{k_p},
          \end{eqnarray*}        
 where  $    K \equiv k_1 + \cdots + k_p$, and
          \begin{eqnarray*}
          ({\rm II}) ~~~~        c_p^{(X)} &=&  a_0^{(p+1)X} \left|      
          \left(  \begin{array}{cccc}
          a_0 & ~ & ~ & ~ \\
          a_1& a_0  &~ \mbox{{\rm \Huge 0}}  & ~  \\
          \vdots & ~& \ddots ~~&~ \\
          a_{p} & a_{p-1}  & \cdots  & a_0  \end{array}   
         \right)^{\mbox{ -X}}+ \left( \left. \begin{array}{cccc}
            ~ & ~  & ~&~  \\
            ~ &  \mbox{{\rm \Huge 0}} & ~ &~ \\
            ~ & ~ & ~ &~\\
            ~ &  ~ & ~&~  \end{array} \right|
         \begin{array}{c} 
             ~1~\\
             ~0~\\
             \vdots\\
            - 1/a_0^X  \end{array} \right)    \right|.
          \end{eqnarray*}   
          \end{I}  
 The proof of (II) requires some additional machinery and will be postponed to Section II.\\ 
$Proof~ of ~(\rm {I})$:  Since $(  a_0 + a_1z + \cdots ) ^X-  (  a_0 + a_1z + \cdots + a_pz^p)^X \sim z^{p+1}$, 
the coefficient of the $z^p$ term is the same in the two sums, so in determining the $c_p^{(X)}$ coefficient we need deal only 
with the finite sum.  Using the generalized binomial theorem, we write         
          \begin{eqnarray}
          \left(  \sum_{q=0}^{p} a_qz^q \right) ^X 
          &=&  \sum_{N=0}^{\infty} \left( \begin{array} {c} X \\ N \end{array} \right)  a_0^{X-N} \left(  \sum_{q=1}^{p}
          a_qz^q \right) ^N \nonumber\\
          &=&   \sum_{N=0}^{\infty} \left( \begin{array} {c} X \\ N \end{array} \right)  a_0^{X-N} 
          \sum_{ 0 \leq k_1, \ldots, k_p \leq N}   \left( \begin{array}{c} N \\ k_1, \ldots, k_p \end{array} \right)
          a_1^{k_1} \cdots a_p^{k_p}~z^{\sum mk_m} .
          \end{eqnarray} 
In the second line we've applied a multinomial expansion to $(\sum_{q=1}^p a_qz^q)^N$.
We now write  $z^{\sum mk_m}  = \sum_{s=0}^{\infty} z^s \delta_{s, \sum mk_m}$ and interchange the order of the 
 sums:
           \begin{eqnarray}
           \left(  \sum_{q=0}^{p} a_qz^q \right) ^X  
           & =&  \sum_{s=0}^{\infty} z^s  \sum_{ 0 \leq k_1, \ldots, k_p \leq s} \sum_{N=0}^{\infty} 
           \left( \begin{array} {c} X \\ N \end{array} \right)
           \left( \begin{array}{c} N \\ k_1, \ldots, k_p \end{array} \right)  
           \delta_{s, \sum mk_m} ~a_0^{X-N}  a_1^{k_1} \cdots a_p^{k_p} 
           \end{eqnarray}
 After the interchange, the $k_q$'s are fixed in the sum over $N$ and only the $N =  k_1+\cdots + k_p~( \equiv K)$
  term contributes.    We are interested in the coefficient of $z^p$ in this sum, which is thus
            \begin{eqnarray}
            c_p^{(X)} &=&    \sum_{ 0 \leq k_1, \ldots, k_p \leq p} 
            \left( \begin{array} {c} X \\ K \end{array} \right)
           \left( \begin{array}{c} K \\ k_1, \ldots, k_p \end{array} \right)  
           \delta_{p, \sum mk_m} ~a_0^{X-K}  a_1^{k_1} \cdots a_p^{k_p}. 
           \end{eqnarray}        
As stated above, this is also the coefficient of $z^p$ in the infinite sum.  QED.
 
 The Kroneker delta in eq.(6) restricts the sums over the $k$'s to a sum over the partitions of $p$.  The multinomial coefficient 
 in this formula counts the number of unique ways the parts can be ordered.  As a sum instead over compositions,
 eq.(6) becomes
          \begin{eqnarray} 
          c_p^{(X)} &=&   \sum_{K=0}^p  \left( \begin{array} {c} X \\ K \end{array} \right)
          \sum_{ 1 \leq q_1, \ldots, q_K \leq p}    
          \delta_{p, q_1+ \cdots +q_K} ~a_0^{X-K}  a_{q_1} \cdots a_{q_K} ;
          \end{eqnarray}
As a corollary, we have the multinomial identity:           
          \begin{4} For fixed $p$,
          \begin{eqnarray*}
  	 \sum_{K=0}^p (-1)^K \sum_{1 \leq q_1, \ldots , q_K \leq p}   \left( \begin{array}{c} p \\
           q_1, \ldots , q_K \end{array} \right) =(-1)^p
          \end{eqnarray*}           
          \end{4}
$Proof$: We set $X=-1$ and $a_q =1/q!$, so that $(\sum a_qz^q)^X = e^{-z}$ and  so $ c_p^{(-1)}  = (-1)^p/p!$.
Then, from eq, (7),
           \begin{eqnarray}
           \frac{(-1)^p}{p!} &=&   \sum_{K=0}^p  (-1)^K \sum_{ 1 \leq q_1, \ldots, q_K \leq p}
           \frac{ \delta_{p, q_1+ \cdots +q_K}}{q_1! \cdots q_K!}  .
           \end{eqnarray} 
Multiplying both sides by $p!$ gives the result.  \\
QED

We now apply Part (I) of Theorem 1 to eq.(1) and use Euler's pentagonal theorem:
          \begin{eqnarray}
          \prod_{k=1}^{\infty} (1-z^k) =   1-z -z^2 + z^5 + z^7 -z^{12} -\cdots .
	 \end{eqnarray}
where the exponents $0,1,2,5,7,12, \ldots $ are generalized pentagonal numbers (sequence A001318): 
$q_m = (3m^2-m)/2,~ m = 0,\pm 1, \pm 2, \cdots$  \cite{Partition}.  $p(n)$ is then equal to a sum over 
the pentagonal partitions of $n$:         
          \begin{eqnarray}
          p(n) &=&    \sum_{ 0 \leq k_1, k_2, k_5,\ldots  k_{q_{M}}  \leq n} (-1)^{\sum k_{q_{2m}} } 
          \left( \begin{array}{c} K \\ k_1, k_2, k_5,  \ldots ,k_{q_M} \end{array} \right) 
          \delta_{n,  \sum k_{q_{m}} q_m}  
          \end{eqnarray}
 where $q_M$ is the largest GPN $\leq n$.   Eq.(10) thus expresses $p(n)$ as a finite sum of integers.  The number of terms in the sum is the number of partitions of $n$ into generalized pentagonal numbers, (sequence  A095699).   For example,
 9 has 10 pentagonal partitions, ($9 = 7+2 = 7+1+1 =\cdots = 1 + \cdots +1$), and $p(9)$ is
         \begin{eqnarray*}
         p(9)  &=&-  \left( \begin{array}{c} 2 \\ 0,1,0,1 \end{array} \right) -  \left( \begin{array}{c} 3 \\ 2,0,0,1 \end{array} \right)
         -  \left( \begin{array}{c} 3 \\ 0,2,1,0 \end{array} \right) 
         -  \left( \begin{array}{c} 4 \\ 2,1,1,0 \end{array} \right)  -  \left( \begin{array}{c} 5 \\ 4,0,1,0 \end{array} \right)\\
         && + \left( \begin{array}{c} 5 \\ 1,4,0,0 \end{array} \right)  
          + \left( \begin{array}{c} 6 \\ 3,3,0,0 \end{array} \right)  +
          \left( \begin{array}{c} 7 \\ 5,2,0,0 \end{array} \right) 
         +  \left( \begin{array}{c} 8 \\ 7,1,0,0 \end{array} \right) +
         \left( \begin{array}{c} 9 \\ 9,0,0,0 \end{array} \right) .
           \end{eqnarray*}  
  
 From the as-yet-to-be-proven Part II of the theorem,  $p(n)$ is also expressible as a $(k+1) \times (k+1)$ determinant, 
 for any integer $k \geq n$;
        \begin{eqnarray}
       p(n) &=& \left| \begin{array} {ccccccc}
       ~~1 & ~~0 & ~~0 & ~~0& ~~\cdots &  ~~0 &~~0 \\
       -1 & ~~1 & ~~0 & ~~0 &~&~& ~~\vdots \\
        -1  & -1 & ~~1 &~~0 &~& ~ & ~~0 \\
       ~~0  & -1 & -1 & ~~1&~& ~& ~~1 \\
        ~~\vdots  & ~ & ~ &~&\ddots & ~& ~~0 \\ 
         d_{k-1}& d_{k-2} & d_{k-3} &d_{k-4}&~& ~~1 &~~ \vdots\\
        d_k & d_{k-1} &d_{k-2} &d_{k-3}& \cdots & -1 & ~~0
          \end{array} \right|       
             \begin{array}{c}~ \\ ~ \\  ~\\ ~ \\
             \left. \begin{array}{c} ~ \\  ~\\ ~ \end{array} \right\} \\
             ~   \end{array}
             \begin{array}{c} ~\\~\\~ \\~\\ n\\ ~   \end{array}            
          \end{eqnarray} 
   where, for $0 \leq q \leq k,$ 
         \begin{eqnarray*}
         d_q = ({\rm sequence~A010815}) = \left\{ \begin{array} {l} (-1)^m ~{\rm if~} q=q_m, ~m= 0, \pm 1, \pm 2, \ldots , \\
             ~~~ 0 ~~{\rm otherwise}, \end{array} \right.  ,             
          \end{eqnarray*}  
This reduces to the form stated in the theorem by successive expansions  by minors along the top row.    Form (11)  will prove more useful in some of the following discussions, but it can be reduced (by expansions in minors) to the $n \times n$  determinant,         
    \begin{eqnarray} 
    \begin{array}{c}  {\rm GPN's}   \\  \hline  0\\1\\2\\~\\~\\5\\~\\7 \\~\\~\\~\\~\\12\\  \vdots \\~\\~\\~
    \end{array} ~~~~
  p(n) = \left| ~~\begin{array}{cccccccccccccccc}
   ~1& -1 & ~\\
   ~1& ~1& -1 & ~\\
   ~0& ~1& ~1& -1 & ~\\
    ~0& ~0& ~1& ~1& -1 & ~\\ 
     -1 & ~0& ~0& ~1&~1 & -1 & ~\\
     ~0& -1 & ~0& ~0&  ~ 1&~1&  -1 & ~\\
       -1&  ~0& -1 & ~0& ~0&  ~ 1&~1&  -1 & ~\\
       ~0& -1& ~0&  -1& ~0 &~ 0&   ~1&~1&  -1 & ~\\
      ~0& ~0& -1& ~0& -1 & ~0& ~0&  ~ 1&~1&  -1 & ~\\ 
       ~0& ~0& ~0& -1 & ~ 0& -1 & ~0& ~0& ~  1&~1&  -1 & ~\\
        ~  0& ~0&~ 0& ~0& -1 & ~0& -1 & ~0& ~0&   ~1&~1&  -1 & ~\\
      ~1 & ~0&~ 0& ~0& ~ 0& -1 &  ~0& -1 & ~0& ~0&  ~ 1&~1&  -1&~ \\
        ~0& ~1&~ 0& ~0&~ 0& ~0& -1 & ~ 0& -1 &~ 0&~ 0&  ~ 1&~1&  ~\\
         ~  \vdots & ~ & ~&~&~&~&~&~&~&~&~&~&~&~&\ddots &        
    \end{array} ~~\right| _{(n \times n)}
       \end{eqnarray}       
      
 \section{Matrix formalism}
 Matrices which are constant along all diagonals are Toeplitz matrices.  We will be concerned in this 
section with lower-triangular Toeplitz (LTT) matrices.   A $nondegenerate$ LTT matrix is one with 
non-zero diagonal elements;  if the diagonal elements are equal to 1, it is then a $unit$ LTT matrix.

Infinite-dimensional LTT matrices have the form 
            \begin{eqnarray}
            A = \left( \begin{array}{ccccc}
            a_0 & ~&~&~&~ \\
            a_1 & a_0 & ~& ~~\mbox{ \Huge 0} &~\\
            a_2 & a_1 & a_0 &~&~  \\
            \vdots & ~ & ~& ~& \ddots
            \end{array} \right) .
            \end{eqnarray}   
The determinant of a $k$-dimensional lower-triangular matrix is the product of it's diagonal elements; if the 
matrix is also Toeplitz, then its determinant is $a_0^k$.
                   
We define the infinite-dimensional lower shift matrix  $J$ as
           \begin{eqnarray}
            J  \equiv \left( \begin{array}{ccccc}
           ~0 ~& ~&~&~&~ \\
           ~1 ~& ~0~ & ~&  \mbox{\Huge 0} &~\\
           ~0 ~& ~1~ &~ 0 &~&~ \\
           \vdots & ~ & ~& ~& \ddots
           \end{array} \right).
            \end{eqnarray}  
The $p$th power of $J$ has elements $(J^p)_{ij} = \delta_{p,i-j} $; these matrices obey the relations
 $J^pJ^q = J^qJ^p = J^{p+q}$.  Any infinite-dimensional LTT matrix can be expanded out in non-negative powers of $J$:
            \begin{eqnarray}
            \left( \begin{array}{ccccc}
            a_0 & ~&~&~&~ \\
            a_1 & a_0 & ~& ~~\mbox{ \Huge 0} &~\\
            a_2 & a_1 & a_0 &~&~  \\
            \vdots & ~ & ~& ~& \ddots
            \end{array} \right) = \sum_{q=0}^{\infty} a_q  J^q. 
            \end{eqnarray}  
The product of two infinite-dimensional LTT matrices is                  
       \begin{eqnarray}
          AB =   \sum_{q=0}^{\infty} a_q J^q~ \sum_{s=0}^{\infty} b_s  J^s 
            = \sum_{p=0}^{\infty} J^p \sum _{q=0}^p a_q b_{p-q} = B A
           \end{eqnarray} 
A similar expression can be written for finite-dimensional LTT matrices, with the $J$'s replaced by finite lower shift
matrices.  LTT matrices, whether finite- or infinite-dimensional, thus commute with one another.

    \begin{1}
Let $A$ be an infinite-dimensional  nondegenerate LTT matrix: 
            \begin{eqnarray*}
            A = \sum_{q=0}^{\infty} a_qJ^q , ~~ a_0 \neq 0.
            \end{eqnarray*}          
Then the inverse of $A$ is an LTT matrix with coefficients
            \begin{eqnarray*}
       b_p &=&  \frac{1}{a_0^{k+1}} \left| \begin{array}{cccccc}
          a_0 & 0  & \cdots & ~& 0 &~  \\
          a_1& a_0  & ~ & ~&~  & ~\\
          \vdots & ~& \ddots &~& ~ & ~ \\
          a_{k-1} &~&~&~& a_0&  ~ \\
          a_{k} & a_{k-1}  & \cdots &~&  a_1&~  \end{array} 
            \begin{array}{c} 
            ~0~\\
             \vdots\\
             ~1~ \\
             \vdots\\
            ~0~ \end{array} \right|  \begin{array}{c}
            ~\\
            ~\\
            ~ \\
             \left. \begin{array}{c} ~ \\  ~ \end{array} \right\} \\
             ~   \end{array}
             \begin{array}{c}
             ~\\
             ~\\
             ~ \\
             p\\
             ~   \end{array}                  
          \end{eqnarray*} 
for any integer  $k \geq p$.  \end{1} 
$Proof$:  We have  
                    \begin{eqnarray}
           \sum_{q=0}^{\infty} a_q J^q~ \sum_{s=0}^{\infty} b_s  J^s 
            = \sum_{p=0}^{\infty} J^p \sum _{q=0}^p a_q b_{p-q} .
                       \end{eqnarray} 
But                       
      \begin{eqnarray}
      \sum_{q=0}^p a_{q}b_{p-q} =  \frac{1}{a_0^{k+1}} \left| \begin{array}{cccccc}
          a_0 & 0  & \cdots & ~& 0 &~  \\
          a_1& a_0  & ~ & ~&~  & ~\\
          \vdots & ~& \ddots &~& ~ & ~ \\
          a_{k-1} &~&~&~& a_0&  ~ \\
          a_{k} & a_{k-1}  & \cdots &~&  a_1&~  \end{array} 
            \begin{array}{c} 
            ~0~\\
             \vdots\\
             ~a_0~ \\
             \vdots\\
            ~a_p~ \end{array} \right|    = \left\{ \begin{array} {l} 1 ~~~~{\rm if~} p =0,\\ 0 ~~~~{\rm if~} p>0, \end{array} \right. 
                      \end{eqnarray}
and therefore $\sum_{s=0}^{\infty} b_s  J^s = A^{-1}$.              QED
  
  It is clear that  this proof still holds if the $J$'s are replaced by finite-dimensional lower shift matrices and the sums over $q$ and $s$ are finite.  Thus, the coefficients of the inverse of a finite $k$-dimensional (nondegenerate) LTT matrix are given by 
 the same formula, and we have the result:
  
  \begin{2}  The inverse of a $k$-dimensional LTT matrix is equal to the $k$-dimensional  truncation of the inverse of the corresponding infinite-dimensional LTT matrix.
 \end{2} 
 I.e.,
                 \begin{eqnarray*}
            \left( \begin{array}{ccccc}
            b_0 & ~&~&~&~ \\
            b_1 & b_0 & ~& ~~ \mbox{{\rm \Huge 0}} &~\\
            b_2 & b_1 & b_0 \\
            \vdots &  ~& ~& \ddots  & ~\\
              \end{array} \right)  = \left( \begin{array}{ccccc}
            a_0 & ~&~&~&~ \\
            a_1 & a_0 & ~& ~~\mbox{{\rm \Huge 0}} &~\\
            a_2 & a_1 & a_0 \\
            \vdots & ~ & ~&  \ddots & ~\\
            \end{array} \right) ^{-\mbox {1}}           
            \end{eqnarray*} 
iff 
           \begin{eqnarray*}
            \left( \begin{array}{ccccc}
            b_0 & ~&~&~&~ \\
             b_1 & b_0 & ~& ~~ \mbox{{\rm \Huge 0}} &~\\
            \vdots &  ~& \ddots  &~ & ~\\
              b_{k-1} & b_{k-2} & \cdots & ~& b_0
              \end{array} \right)  = \left( \begin{array}{ccccc}
            a_0 & ~&~&~&~ \\
            a_1 & a_0 & ~& ~~\mbox{{\rm \Huge 0}} &~\\
            \vdots & ~ &  \ddots &~ & ~\\
            a_{k-1} & a_{k-2} & \cdots & ~& a_0
            \end{array} \right) ^{-\mbox {1}}           
            \end{eqnarray*}           
 for all $k$.

 \begin{3} If $A= \sum_{q=0}^{\infty} a_pJ^p$ is an infinite-dimensional nondegenerate LTT matrix and $X$ is a complex number, then $A$ raised to the power $X$ is
         \begin{eqnarray*}
      A^X =  \sum_{s=0}^{\infty} J^s  \sum_{ 0 \leq k_1, \ldots, k_s \leq s}
       \left( \begin{array} {c} X \\ K \end{array} \right)
           \left( \begin{array}{c} K \\ k_1, \ldots, k_s \end{array} \right)  
            \delta_{s, \sum mk_m} ~a_0^{X-K}  a_1^{k_1} \cdots a_s^{k_s} 
    \end{eqnarray*} 
     \end{3} $Proof$:  As before, we first consider the finite-sum case.  Let  $z$ be some nonzero complex number.  
     We then define ${\bf Z }$ as the infinite row vector
            \begin{eqnarray}
            {\bf Z } = (1, ~z,~z^2,~z^3, \ldots ). 
            \end{eqnarray}                        
Multiplying  {\bf Z } by $J$ on the right, we have  ${\bf Z} J = z {\bf Z}$.   Therefore, ${\bf Z} J^q = z^q {\bf Z}$ and, by linearity, $ {\bf Z} (\sum c_q J^q) =( \sum  c_q z^q) {\bf Z}$.   Then, if $F(z)$ is a holomorphic, nonzero function of $z$ in a neighborhood $U$ of
the origin, $ {\bf Z} F(J) = F(z) {\bf Z}$ for all $ z \in U$, where $F(J)$ is the LTT matrix obtained by replacing $z$ by $J$
in the Taylor series expansion of $F(z)$.  Since $ (\sum_{q=0}^{p} a_pz^p)^X ,~ a_0 \neq 0$ is such a function, we can write, 
     \begin{eqnarray}
      {\bf Z} \left( \sum_{q=0}^p a_qJ^q \right)^X  &=& \left( \sum_{q=0}^{p} a_qz^q \right)^X {\bf Z}  \\
     &=&   \sum_{s=0}^{\infty} z^s  \sum_{ 0 \leq k_1, \ldots, k_p \leq s} 
       \left( \begin{array} {c} X \\ K \end{array} \right)
           \left( \begin{array}{c} K \\ k_1, \ldots, k_p \end{array} \right) 
            \delta_{s, \sum mk_m} ~a_0^{X-K}  a_1^{k_1} \cdots a_p^{k_p} ~  {\bf Z} \nonumber \\
     &=&  {\bf Z}  ~ \sum_{s=0}^{\infty} J^s  \sum_{ 0 \leq k_1, \ldots, k_p \leq s} 
       \left( \begin{array} {c} X \\ K \end{array} \right)
           \left( \begin{array}{c} K \\ k_1, \ldots, k_p \end{array} \right)  
            \delta_{s, \sum mk_m} ~a_0^{X-K}  a_1^{k_1} \cdots a_p^{k_p}  .\nonumber 
      \end{eqnarray}
$z$ is otherwise arbitrary, so we have
     \begin{eqnarray}
      \left( \sum_{q=0}^p a_qJ^q \right)^X    &=&    \sum_{s=0}^{\infty} J^s  \sum_{ 0 \leq k_1, \ldots, k_p \leq s} 
       \left( \begin{array} {c} X \\ K \end{array} \right)
           \left( \begin{array}{c} K \\ k_1, \ldots, k_p \end{array} \right) 
            \delta_{s, \sum mk_m} ~a_0^{X-K}  a_1^{k_1} \cdots a_p^{k_p} .
      \end{eqnarray}
 For fixed $s < p$, $k_{s+1} = \cdots =k_p =0$  in the sums over the $k$'s as a result of the restriction $ s=k_1+ \cdots + pk_p$:
     \begin{eqnarray}
      \left( \sum_{q=0}^p a_qJ^q \right)^X    &=&    \sum_{s=0}^{p-1} J^s  \sum_{ 0 \leq k_1, \ldots, k_s \leq s} 
       \left( \begin{array} {c} X \\ K \end{array} \right) \left( \begin{array}{c} K \\ k_1, \ldots, k_s \end{array} \right)  
            \delta_{s, \sum mk_m} ~a_0^{X-K}  a_1^{k_1} \cdots a_s^{k_s}  \nonumber \\
          &+&    \sum_{s=p}^{\infty} J^s  \sum_{ 0 \leq k_1, \ldots, k_p \leq s} 
       \left( \begin{array} {c} X \\ K \end{array} \right) \left( \begin{array}{c} K \\ k_1, \ldots, k_p \end{array} \right)
            \delta_{s, \sum mk_m} ~a_0^{X-K}  a_1^{k_1} \cdots a_p^{k_p} . 
                \end{eqnarray}
The equality in the lemma is then demonstrated by letting $p \rightarrow  \infty$ in this equation.  QED\\ 
 ~\\
 We now have the necessary machinery to prove the second part of Theorem 1.\\
~\\ 
$Proof~ of ~Thm~1, {\rm (II)}$: 
 From the proof of Lemma 3, the coefficients in the expansion of  $ \left(\sum a_p J^p \right)^X$ in powers of $J$ are the same as in the expansion of $  \left(\sum  a_p z^p \right)^X$ in powers of $z$; i.e.,
      \begin{eqnarray}
      \sum_{q=0}^{\infty} c_q z^q  =   \left( \sum_{q=0}^{\infty} a_q z^q \right)^X  ~~
      {\rm iff}~~  \sum_{q=0}^{\infty} c_q J^q  =   \left( \sum_{q=0}^{\infty} a_q J^q \right)^X 
           \end{eqnarray}          
We can therefore use Lemma 1, written in the form         
         \begin{eqnarray}
       b_p &=&  \frac{1}{a_0^{p+1}} \left| \left( \begin{array}{cccccc}
          a_0 & ~  & ~ & ~& ~ &~  \\
          a_1& a_0  & ~ &  \mbox{{\rm \Huge 0}}   &~  & ~\\
          \vdots & ~& \ddots &~& ~ & ~ \\
          a_{p-1} &~&~&~& a_0&  ~ \\
          a_{p} & a_{p-1}  & \cdots &~&  a_1& ~a_0 \end{array} 
            \right) 
            + \left( \left. \begin{array}{cccc}
            ~ & ~  & ~&~  \\
            ~ & ~  & ~&~  \\
            ~ &  \mbox{{\rm \Huge 0}} & ~ &~ \\
            ~ & ~ & ~ &~\\
            ~ &  ~ & ~&~  \end{array} \right|
             \begin{array}{c} 
             ~1~\\
             ~0 ~\\
              ~\vdots ~\\
              0 \\
              -a_0 \end{array} \right)      \right| ,              
          \end{eqnarray}      
(taking  $k=p$), with the replacement $ A= \sum a_p J^p \rightarrow A= \left(\sum a_p J^p \right)^{-X}$
(and thus $a_0 \rightarrow 1/a_0^X$), to solve for $c_p^{(X)}$ as a determinant in the form (II) in the theorem.  QED.
     
In a straightforward fashion, eq. (23) generalizes to
          \begin{eqnarray}
          \sum_{p=0}^{\infty} c_p z^p &=&  \left(  \sum_{p=0}^{\infty} a_p^{(1)}z^p \right) ^{X_1}
          \cdots  \left(  \sum_{p=0}^{\infty} a_p^{(n)}z^p \right) ^{X_n} \nonumber \\           
        {  \rm iff}~~\sum_{p=0}^{\infty} c_p J^p &=&  \left(  \sum_{p=0}^{\infty} a_p^{(1)}J^p \right) ^{X_1}
          \cdots  \left(  \sum_{p=0}^{\infty} a_p^{(n)}J^p \right) ^{X_n}   . 
          \end{eqnarray} 
The $c_p$ coefficients are then:               
          \begin{eqnarray}
          c_p &=& \left( (a_0^{(1)})^{X_1} \cdots (a_0^{(n)})^{X_n} \right)^{p+1}  \left|    
          \left(  \begin{array}{cccc}
          a_0^{(1)} & ~ & ~ & ~ \\
          a_1^{(1)}& a_0^{(1)}  &~ \mbox{{\rm \Huge 0}}  & ~  \\
         \vdots & ~& \ddots ~~&~ \\
          a_{p}^{(1)} & a_{p-1}^{(1)}  & \cdots  & a_0^{(1)}  \end{array}   
         \right)^{-X_1}\cdots 
         \left(  \begin{array}{cccc}
          a_0^{(n)} & ~ & ~ & ~ \\
          a_1^{(n)}& a_0^{(n)}  &~ \mbox{{\rm \Huge 0}}  & ~  \\
         \vdots & ~& \ddots ~~&~ \\
          a_{p}^{(n)} & a_{p-1}^{(n)}& \cdots  & a_0^{(n)}  \end{array}   
         \right)^{-X_n}\right. \nonumber\\
          && ~~~~~~~~~~~~~~~~~~~~~~~~~~~~~~~~~~\left. 
          + \left( \left. \begin{array}{cccc}
          ~ & ~  & ~&~  \\
          ~ &  \mbox{{\rm \Huge 0}} & ~ &~ \\
          ~ & ~ & ~ &~\\
          ~ &  ~ & ~&~  \end{array} \right|
          \begin{array}{c} 
           ~1~\\
           ~0~\\
           \vdots\\
          - (a_0^{(1)})^{-X_1}\cdots (a_0^{(n)})^{-X_n}  \end{array} \right)    \right|.
          \end{eqnarray}
Gradshteyn and Ryzhik  \cite{Grad} give an equivalent expression for these coefficients for the case 
   $n=2,~X_1=1, ~X_2 =-1$.   
      
Returning to the partition function:  MacMahon's recurrence relation \cite{Partition},
         \begin{eqnarray}
              p(n) -p(n-1) -p(n-2) + p(n-5)+p(n-7)  -p(n-12) - p(n-15) + \cdots =0,
         \end{eqnarray}
follows directly from expression (11); setting $k=n$ in that equation, the sum on the right in (27) is equal to the determinant
        \begin{eqnarray*}
        \left| \begin{array} {ccccccc}
       ~~1 & ~~0 & ~~0 & ~~0& ~~\cdots &  ~~0 &~~1 \\
       -1 & ~~1 & ~~0 & ~~0 &~&~& -1 \\
        -1  & -1 & ~~1 &~~0 &~& ~ & -1 \\
       ~~0  & -1 & -1 & ~~1&~& ~& ~~0 \\
        ~~\vdots  & ~ & ~ &~&\ddots & ~& ~~\vdots \\ 
         ~ & ~ & ~ &~&~& ~~1 &~~d_{n-1}\\
        d_n & d_{n-1} &d_{n-2} &d_{n-3}& \cdots & -1 & ~~d_n
          \end{array} \right| ,                
          \end{eqnarray*} 
  which is zero since the first and last columns are equal.  
  
  Expressed in terms of the $J$  matrices, expression (12) for $p(n)$ is 
\begin{eqnarray}
p(n) &=& \det \left[ -J^T  + \sum_{m>0}^{q_m < n+1} (-1)^{m+1} \left( ~J^{(m-1)(3m+2)/2} + J^{(m+1)(3m-2)/2}~ \right) \right]_{n \times n}
\end{eqnarray}
where $J^T$ is the transpose of $J$ and the notation $[ ~~]_{n \times n}$ means the $n \times n$ truncation of an
infinite-dimensional matrix.  And, by relation (23),  we have the compact matrix equivalent of $p(n)$'s generating-function
equation:         
       \begin{eqnarray}
      \left( \begin{array} {ccccc}
       p(0) & 0 & 0 & \cdots  &~~0 \\
       p(1) & p(0) & 0  &~& ~~0 \\
       p(2)  & p(1) & p(0)  &~&  ~~0 \\
             ~~\vdots  & ~ & ~ &\ddots & ~~\vdots \\ 
        p(k) & p(k-1) & p(k-2) & \cdots  & ~~p(0)
          \end{array} \right) &=& \left( \begin{array} {ccccc}
       ~~1 & ~~0 & ~~0 & ~~\cdots  &~~0 \\
       -1 & ~~1 & ~~0  &~& ~~0 \\
        -1  & -1 & ~~1  &~&  ~~0 \\
        ~~\vdots  & ~  &~~ &\ddots & ~~\vdots \\ 
        d_k & d_{k-1} &d_{k-2} & \cdots  & ~~1
          \end{array} \right)^{-1}      .        
          \end{eqnarray} 
          
 The generating function for the number of partitions in which no part occurs more than $D$ times is \cite {Partition}
\begin{eqnarray}
\prod_{k=1}^{\infty} \frac { 1-x^{(D+1)k} } {1-x^k} = \frac { 1-x^{D+1} - x^{2(D+1)} +  x^{5(D+1)} +  x^{7(D+1)} - \cdots }
{1-x-x^2 + x^5 +x^7 - \cdots }.
\end{eqnarray}
We have then from (26), with $X_1 =-X_2 =1$, 
          \begin{eqnarray}
       p_D(n)        = \left| \begin{array} {ccccc}
       ~~1 & ~~0  & ~~\cdots &  ~~0 &~~ t_0 \\
       -1 & ~~1   &~&~~0& ~~  t_1 \\
          ~~\vdots  & ~  &\ddots & ~& ~~\vdots \\ 
         d_{n-1} & d_{n-2} & ~ & ~~1 &~~ t_{n-1} \\
        d_n & d_{n-1} &  \cdots & -1 & ~~ t_{n}
        \end{array} \right|     , 
            \end{eqnarray}  
 where         
              \begin{eqnarray}
          t_q = \left\{ \begin{array} {l} (-1)^m ~{\rm if~} q=(D+1)q_m, ~m=  0, \pm 1, \pm 2, \ldots  
             ~~~~~~~\\~~~ 0 ~~{\rm otherwise}. \end{array} \right.                
          \end{eqnarray} 
 By expanding this determinant by minors along the last column, and using the expression (11), we get the relation            
           \begin{eqnarray}
          p_D(n) = p(n)  + \sum_{m} (-1)^m p(n- (D+1)q_m) ,
          \end{eqnarray}
 which is a generalization of (27), the  $D=0$ case.  If we now take $D=1$, $p_1(n) \equiv q(n)$ = the number of partitions of $n$ into distinct integers  (A000009).  This is also, from a result due to Euler \cite{PartitionQ},  the number of partitions of $n$ into $odd$ integers, and so we have for the odd partition function
           \begin{eqnarray}
          q(n) = p(n) - p(n-2)-p(n-4) + p(n-10) + \cdots +(-1)^m p(n-2p_m) + \cdots .
          \end{eqnarray}
 		
 \section{Application to Bernoulli, Euler  and Stirling numbers}    
 Vella's expression for the $n$th Bernoulli number,  (eq.(a) in his Theorem 11) is, in my notation,
          \begin{eqnarray}
          B_n 	 &=& \sum_{K=0}^n \frac{ (-1)^K }{1+K}  \sum_{1 \leq q_1, \ldots , q_K \leq n} \left( \begin{array}{c} n \\
	q_1, \ldots , q_K \end{array} \right) .
	\end{eqnarray}
  Applying Part I of Theorem 1 to the generating function (2a) for Bernoulli numbers, we have
           \begin{subequations}
           \begin{eqnarray}
          B_n &=& n!\sum_{0 \leq k_1,  \cdots ,k_n \leq n }
	 \left( \begin{array}{c} K \\ k_1, \ldots , k_n \end{array} \right)  \delta_{n,\sum_{m=1}^n mk_m } 
	 \left( \frac {-1~}{2!} \right)^{k_1}  \left( \frac{-1~}{3!} \right)^{k_2} \cdots  \left( \frac {-1~}{(n+1)!} \right)^{k_n} .
	 \end{eqnarray}
As a sum over compositions, this expression becomes
	 \begin{eqnarray}
	 B_n &=& \sum_{K=0}^n (-1)^K  \sum_{1 \leq q_1, \ldots , q_K \leq n} 
           \frac{1}{(q_1+1) \cdots  (q_K+1)} \ \left( \begin{array}{c} n \\
	q_1, \ldots , q_K \end{array} \right),
          \end{eqnarray}
          \end{subequations}
to be  compared to Vella's result above.  

 Another expression for the (even-numbered) Bernoulli numbers is obtained from the generating function
       \begin{eqnarray}
           \sum _{q=0}^{\infty} \frac{(2-2^{2q})B_{2q}}{(2q)!}~ z^{2q} &=&  \frac{z}{\sinh z}.         
              \end{eqnarray}
From this we have
          \begin{subequations}
	 \begin{eqnarray}
	 B_{2p} &=& \frac{(2p)!}{2-2^{2p}}\sum_{0 \leq k_1, \ldots, k_p \leq p} 
	  \left( \begin{array}{c} K \\ k_1, \ldots , k_p \end{array} \right)
	 \delta_{p,\sum_{m=1}^p mk_m } \left( \frac{-1~}{3!} \right)^{k_1}  \nonumber \\
	 && ~~~~~~~~~~~~~~~~~~~~~~~~~~~~~~~~~~~~~~~~~~~~~~~~\times \left( \frac{-1~}{5!} \right)^{k_2}
		\cdots \left( \frac{-1~}{(2p+1)!} \right)^{k_p} ,  \\
	&=& \frac{1}{2-2^{2p}} \sum_{K=0}^p (-1)^K \sum_{1 \leq q_1, \cdots , q_K \leq p}
	\frac{1}{(2q_1+1) \cdots (2q_K+1)} \left( \begin{array}{c} 2p \\
	2q_1, \ldots , 2q_K \end{array} \right),
	\end{eqnarray}	
	\end{subequations}      
while from the generating function for the Euler numbers (A000364), eq.(2b), we obtain the two expressions,
           \begin{subequations}
	 \begin{eqnarray}
	 E_{2p}&=& (2p)! \sum_{0 \leq k_1, \ldots, k_p \leq p}~  \left( \begin{array}{c} K \\ k_1, \ldots , k_p \end{array} \right)
	\delta_{p,\sum_{m=1}^p mk_m }  \left( \frac{-1~}{2!} \right)^{k_1} \left( \frac{-1~}{4!} \right)^{k_2}
	\cdots \left( \frac{-1~}{(2p)!} \right)^{k_p} , \\
	&=& \sum_{K=0}^p (-1)^K \sum_{1 \leq q_1, \ldots , q_K \leq p} \left( \begin{array}{c} 2p \\
	2q_1, \ldots , 2q_K \end{array} \right).
	\end{eqnarray}
	\end{subequations}

 $B_{2p}$ and $E_{2p}$ can thus both be expressed as sums over the even compositions of $2p$.   Since there are no even 
 compositions of odd numbers, these expressions can be extended to include all of the odd-numbered numbers 
 except for $B_1$.  We then have, more generally,
 	 \begin{eqnarray}
	 B_{n} &=& \frac{1}{2-2^{n}} \sum_{K=0}^{ \lfloor n/2 \rfloor } (-1)^K \sum_{1 \leq q_1, \ldots , q_K \leq \lfloor n/2 \rfloor } 
	 \frac{1}{(2q_1+1) \cdots (2q_K+1)}	 \left( \begin{array}{c} n \\
	2q_1, \ldots , 2q_K \end{array} \right),\\ && ( n \neq 1) ,\nonumber
	\end{eqnarray}	 
	 \begin{eqnarray}
	 E_{n} &=& \sum_{K=0}^{ \lfloor n/2 \rfloor } (-1)^K \sum_{1 \leq q_1, \ldots , q_K \leq \lfloor n/2 \rfloor }
	  \left( \begin{array}{c} n \\ 2q_1, \ldots , 2q_K \end{array} \right)
	\end{eqnarray}	  
 This expression for $E_{n}$, in a different notation, was previously derived by Vella, (eq. (c) in Thm. 11 \cite{Vella}).   
  $E_n$  can however also be expressed as a sum over the $odd$ partitions/compositions of $n-1$: 
        \begin{II} For n $>$1,
            \begin{eqnarray*}
	 E_{n} = \sum_{N=1}^{\lfloor n/2 \rfloor } (-1)^N \sum_{1 \leq q_1, \ldots ,q_{2N-1} \leq \lfloor n/2 \rfloor }
	  \left( \begin{array}{c} n-1\\ 2q_1-1, \ldots , 2q_{2N-1}-1 \end{array} \right)	 
        \end{eqnarray*}
        \end{II}
~\\   
  $Proof.$ The sum is over all compositions of $n-1$ that contain an odd number of odd parts, which is an empty set if $n-1$ is
  even.   Therefore $E_n =0$ for odd $n$, and we only have to prove this equation for even $n$.  The equality is trivially true for $n=2$.  We will show that both $E_{2p}$ and               
         \begin{eqnarray}
	 a_{2p} \equiv \sum_{N=1}^p (-1)^N \sum_{1 \leq q_1, \ldots ,q_{2N-1} \leq p} \left( \begin{array}{c} 2p-1\\ 
	2q_1-1, \ldots , 2q_{2N-1}-1 \end{array} \right) 	 
        \end{eqnarray}
satisfy the same recursion relation  
          \begin{eqnarray}
	a_{2p} = - 1 -\sum_{q=1}^{p-1} \left( \begin{array}{c} 2p-1 \\ 2q \end{array} \right)  2^{2q-1} a_{2p-2q}.
	 \end{eqnarray}
Then the equality would, by induction, be valid for all $p$.  
 
 To prove this for $E_{2p}$, consider the expression 
         \begin{eqnarray*}
          \frac{d}{dz} ~\frac{1}{\cosh z} + 2 \sinh z + \cosh 2z ~ \frac{d}{dz} ~\frac{1}{\cosh z} .
           \end{eqnarray*}
This is equal to zero as a result of the identity $\cosh 2z = 2 \cosh ^2z-1$.  Expanding out $\sinh z$ and $1/\cosh z$ and performing the differentiation in the 1st and 3rd terms, we have, for each power of $z$,
	 \begin{eqnarray}
	\frac{ E_{2p}}{(2p-1)!} + \frac{ 2}{(2p-1)!} + \sum_{q=0}^{p-1}  \frac{2^{2q} E_{2p-2q}}{(2q)!(2p-2q-1)!}  =0,
	 \end{eqnarray}
which can be rearranged into the form of the recursion relation (43) above.

On the other hand  we have
  	 \begin{eqnarray}
	 a_{2p}  &=&- 1 + \sum_{N=2}^p (-1)^N
	\sum_{q_1=1}^{p-1}  \sum_{q_2=1}^{p-1} \sum_{k_1=1}^{p-1} \cdots \sum_{k_{2N-3}=1}^{p-1}
         \left( \begin{array}{c} 2p-1 \\ 2q_1-1, 2q_2-1, 2k_1-1, \ldots , 2k_{2N-3}-1 \end{array} \right), ~\nonumber \\
         \end{eqnarray}
where we've separated off the $N=1$ term and made the substitutions $q_3 \rightarrow k_1, \ldots,
        q_{2N-1} \rightarrow k_{2N-3}$.   We have then, with $N = L+1$,
         \begin{eqnarray} 
	a_{2p} &=& - 1+\sum_{q_1=1}^{p-1}\sum_{q_2=1}^{p-1} \frac{1} {(2q_1-1)!(2q_2-1)!}~ \frac{  (2p-1)! }
	{(2p-2q_1-2q_2+1)!} \nonumber\\
	&& ~~~~~~~~~ \times \sum_{L=1}^{p-q_1-q_2 +1} (-1)^{L+1} \sum_{k_1=1}^{p-q_1-q_2+1} 
	\cdots \sum_{k_{2L-1}=1}^{p-q_1-q_2+1} 
	\left( \begin{array}{c} 2p-2q_1-2q_2+1 \\2k_1-1, \ldots , 2k_{2L-1}-1 \end{array} \right).
	    \end{eqnarray}	    
Setting $q =q_1+q_2-1$, this is     
         	 \begin{eqnarray}
         a_{2p} &=& - 1 -\sum_{q=1}^{p-1} ~\sum_{q_1=1}^{q}  \left( \begin{array}{c} 2q \\ 2q_1-1 \end{array} \right) ~ 
         \frac{(2p-1)!}{(2q)!(2p-2q-1)! } \nonumber \\
	 && ~~~~~~~~~~~~~ \times \sum_{L=1}^{p-q} (-1)^{L} \sum_{k_1=1}^{p-q} 
	\cdots \sum_{k_{2L-1}=1}^{p-q} \left( \begin{array}{c} 2p-2q-1 \\2k_1-1, \ldots , 2k_{2L-1}-1 \end{array} 
	\right).
	 \end{eqnarray}	    
The sum over $q_1$ equals $2^{2q-1}$, evaluated by expanding $(1+1)^{2q} -(1-1)^{2q}$ binomially.  
The remaining sum over $L$ and $k_1, \ldots, k_{2L-1}$ is $a_{2p-2q}$, and we again get relation (43).
   QED 
   
The sum-over-partitions form follows in a straightforward fashion:  
            \begin{eqnarray}
	 E_{2p} &=&  (-1)^{p-1} (2p-1)! \sum_{0 \leq k_1, \ldots, k_p \leq 2p-1} 
	 \left( \begin{array} {c} K \\ k_1, \ldots , k_p \end{array} \right)
	\delta_{2p-1,\sum (2m-1)k_m } \nonumber \\
	&& ~~~~~~~~~~~~~~~~~~~~~~~~~~~ \times  \left( \frac{-1~}{1!} \right)^{k_1}  \left( \frac{1}{3!} \right)^{k_2}
	   \cdots \left( \frac{(-1)^p}{(2p-1)!} \right)^{k_p}   .
	 \end{eqnarray}
The sign of each term in this sum is $(-1)^{p-1+ \sum mk_m}$, while the sign in the sum-over-compositions form is
 $(-1)^N$.    To check that the signs agree, note that $2p-1 = \sum (2m-1)k_m = 2 \sum mk_m -K$.  But $K$ is the number of
factorials in the denominator of each term, which in the sum-over-compositions form is $2N-1$.  So  $ \sum mk_m = N+p-1$
and $(-1)^{p-1+ \sum mk_m} = (-1)^N$.
 
As an example of the even and the odd expansions for $E_{n}$,  we have:
	 \begin{eqnarray}
	 E_{10} &=&10! \left( - \frac{1}{10!} + \frac{2}{2!8!} + \frac{2}{4!6!}
	- \frac{3}{2!^2 6!}- \frac{3}{2!4!^2} +\frac{4}{2!^3 4!} - \frac{1}{2!^5}\right) \nonumber \\  
	 &=&~9! \left( - \frac{1}{9!} + \frac{3}{1!^27!} + \frac{6}{1!3!5!}
	+\frac{1}{3!^3}- \frac{5}{1!^45!} -\frac{10}{1!^33!^2} + \frac{7}{1!^6 3!} - \frac{1}{1!^9}\right) = -50,521
		\end{eqnarray}  

In a similar fashion, Theorem 1 can be applied to the Euler and Bernoulli polynomials, using their generating functions:
          \begin{subequations}
          \begin{eqnarray} 
          \sum_{q=0}^{\infty} \frac {E_q(x)}{q!}~ z^q &=& \frac{2e^{xz}}{e^z+1} 
          = 2\left\{  \sum _{k=0}^{\infty} z^k~ \frac{(1-x)^k+(-x)^k}{k!} \right\}^{-1},\\
          \sum_{q=0}^{\infty} \frac{ B_q(x)}{q!} ~z^q &=& \frac{ze^{xz}}{e^z-1} 
          =  \left\{ \sum_{k=0}^{\infty} z^k ~\frac{( 1-x)^{k+1} - (-x)^{k+1}}{(k+1)!} \right\}^{-1}.
            \end{eqnarray}
	  \end{subequations}
The results are straightforward and we omit writing out the explicit expressions.

Bell polynomials \cite{Bell} are defined as
          \begin{eqnarray}
          B_{n,k}(x_1, \ldots , x_{n-k+1})  &=&  \sum_{0\leq k_0, \ldots, k_{n-k} \leq n} \delta_{k,k_0+\cdots +k_{n-k}}  
          \delta _{n,\sum_m    mk_m} \nonumber \\
          && ~~~~~~~\times \frac{n!}{k_0! \cdots k_p!}~\left( \frac{x_1}{1!} \right) ^{k_0} \cdots 
          \left( \frac{x_{n-k+1}}{(n-k+1)!} \right)^{k_{n-k}}  .
         \end{eqnarray}
Stirling numbers of the 2nd kind are equal to the values of these polynomials at $x_1=x_2= \cdots =x_k=1$: 
          \begin{eqnarray}
          S(n,n-p) &=& \left. B_{n,n-p}(x_1,x_2, \ldots , x_{p+1}) \right|_{x_1= \cdots =x_{p+1} =1} \nonumber \\ 
          &=& \frac{n!}{ (n-p)!} \sum_{0\leq k_1, \ldots, k_p \leq p}  \left( \begin{array}{c} n-p\\K \end{array} \right) 
          \left( \begin{array} {c} K \\ k_1, \ldots , k_p \end{array} \right)
          \delta _{p,\sum  mk_m} \nonumber \\
          && ~~~~~~~~~~~~~~~~~~~~~~~~~~~~~~~~~~~~~~~~~~\times   \left( \frac{1}{2!} \right)^{k_1} \left( \frac{1}{3!} \right)^{k_2} 
          \cdots \left( \frac{1}{(p+1)!} \right)^{k_p} , 
          \end{eqnarray}
where we've set $k=n-p$ and in the last line we've summed over $k_0$.   A similar expression for Stirling 
numbers of the 1st kind can be found using their relation to the $n$-th order Bernoulli numbers \cite{Korn},
             \begin{eqnarray}
             s(n,n-p)= \left( \begin{array}{c}
             n-1\\
             p \end{array} \right) B_{p}^{(n)},
             \end{eqnarray}              
which have the generating function
   \begin{eqnarray}
 \sum_{q=0}^{\infty} \frac{B_q^{(k)}}{q!} ~ z^q  =  \frac{z^k}{(e^z-1)^k}  .
   \end{eqnarray}
Using Theorem 1, we get for $s(n,n-p)$:
         \begin{eqnarray}
             s(n,n-p)  &=& \frac{(n-1)!}{(n-p-1)!}\sum_{0 \leq k_1, \ldots , k_{p}  \leq p}  
                           \left( \begin{array}{c} n+K-1\\K \end{array} \right) 
                            \left( \begin{array} {c}  K \\ k_1, \ldots , k_p \end{array} \right)   \delta_{p,\sum  mk_m} \nonumber  \\
             && ~~~~~~~~~~~~~~~~~~~~~~~~~~~~~~ \times  \left( \frac{-1~}{2!} \right)^{k_1} 
             \left( \frac{-1~}{3!} \right)^{k_2}  \cdots  \left( \frac{-1}{(p+1)!} \right)^{k_p} .              
                  \end{eqnarray}  
                  
  Laplace's formula (3) for the Bernoulli numbers corresponds to expression (36a) and follows from equation (2a) and Part (II) of Theorem 1.  Other matrix representations we can derive from this theorem are:                 
             \begin{subequations}
          \begin{eqnarray}
          B_{2p} &=& - \frac{(2p)!}{2^{2p}-2}  \left| \begin{array}{cccccc}
            1 & 0 & 0 & ~ & \cdots &   ~ 1 \\
            \frac{1}{3!} & 1 & 0 &~& ~ &   ~0 \\
            \frac{1}{5!} & \frac{1}{3!} & 1 & ~ & ~& ~\vdots \\
            \vdots & ~ & ~ & \ddots  &~&  ~ \\
            \frac{1}{(2p+1)!}&\frac{1}{(2p-1)!} &\frac{1}{(2p-3)!} &  \cdots & \frac{1}{3!} ~& 0
            \end{array} \right| ; \\
	 E_{2p} &=&  (2p)! ~~\left| \begin{array}{cccccc}
            1 & 0 & 0 & ~ & \cdots &   ~ 1 \\
            \frac{1}{2!} & 1 & 0 &~& ~ &   ~0 \\
            \frac{1}{4!} & \frac{1}{2!} & 1 & ~ & ~& ~\vdots \\
            \vdots & ~ & ~ & \ddots  &~&  ~ \\
            \frac{1}{(2p)!}&\frac{1}{(2p-2)!} &\frac{1}{(2p-4)!} &  \cdots & \frac{1}{2!} ~& 0
            \end{array} \right| ; \\
             S(n,n-p) &=& \frac{n!}{(n-p)!}  ~\left| ~I + \left(  \begin{array}{cccc}
               1 & ~ & ~ & ~ \\
              \frac{1}{2!}& 1 &~~~~ \mbox{{\rm \Huge 0}}  & ~  \\
              \vdots & ~& ~~\ddots &~ \\
              \frac{1}{(p+1)!} &  \frac{1}{p!}   & ~~\cdots  & ~~~~~1  \end{array}          
              \right)^{\mbox{\it{n-p}}}   \times
             \left( \left. \begin{array}{cccc}
            ~ & ~  & ~&~  \\
            ~ &  \mbox{{\rm \Huge 0}} & ~ &~ \\
            ~ & ~ & ~ &~\\
            ~ &  ~ & ~&~  \end{array} \right|
             \begin{array}{c} 
             ~1~\\
              ~0~\\
              \vdots\\
             - 1 \end{array} \right)  ~ \right|;\\ 
              s(n,n-p)  &=&  \frac{(n-1)!}{(n-p-1)!}           
              \left|~ \left( \begin{array}{cccc}
              1 & ~ & ~ & ~ \\
              \frac{1}{2!}& 1 &~~~~ \mbox{{\rm \Huge 0}}  & ~  \\
              \vdots & ~& ~~\ddots &~ \\
              \frac{1}{(p+1)!} &  \frac{1}{p!}   & ~~\cdots  & ~~~~~1  \end{array}             
              \right)^{\mbox{\it{n}}}  +      
             \left( \left. \begin{array}{cccc}
            ~ & ~  & ~&~  \\
            ~ &  \mbox{{\rm \Huge 0}} & ~ &~ \\
            ~ & ~ & ~ &~\\
            ~ &  ~ & ~&~  \end{array} \right|
             \begin{array}{c} 
             ~1~\\
              ~0~\\
              \vdots\\
             - 1 \end{array} \right)  ~  \right|.
            \end{eqnarray}
            \end{subequations}           
 Note that the size of the determinant for $B_{2p}$ in (56a) is $(p+1)\times(p+1)$, (as is the one for $E_{2p}$ 
in (56b)), compared with $(2p+1)\times(2p+1)$ using Laplace's formula.

\section {Conclusion}
We have derived compact, closed-form expressions for partition functions and for Bernoulli, Euler and Stirling numbers
that contain only ``simple'' numbers and that require either finite summations or finding the determinants or inverses of 
matrices.  In particular, the partition function $p(n)$ is given directly by a sum of integers, the number of terms in the sum being the value at $n$ of  the pentagonal partition function.  

 \begin{acknowledgments}
I wish to thank Yonko Millev for invaluable assistance.
\end{acknowledgments}
            
\appendix{}
\section{ $p(5k+4), ~p(7k+5)$, and $p(25k+24)$ determinants}
For $n= 5k+4, 7k+5$, or  $25k+24$, the dimension of the matrix in eq.(12) can be reduced to $k+1$ by using one of the 
identities below by Ramanujan \cite{Partition,Ram}: 
          \begin{subequations}
          \begin{eqnarray}
          \sum_{k=0}^{\infty} p(5k+4)q^k &=& 5~ \frac{(q^5)^5_{\infty}}{(q)^6_{\infty} },\\
          \sum_{k=0}^{\infty} p(7k+5)q^k &=& 7~ \frac{(q^7)^3_{\infty}}{(q)^4_{\infty} } + 49q~ \frac{(q^7)^7_{\infty}}{(q)^8_{\infty} },\\    
          \sum_{k=0}^{\infty} p(25k+24)q^k &=& 5^2 \cdot 63  \frac{(q^5)^6_{\infty}}{(q)^7_{\infty} }  
           + 5^5 \cdot 52q~ \frac{(q^5)^{12}_{\infty}}{(q)^{13}_{\infty} }   +5^7 \cdot 63q^2~ \frac{(q^5)^{18}_{\infty}}{(q)^{19}_{\infty} }
           \nonumber \\
           &&~~~~ + 5^{10} \cdot 6q^3~ \frac{(q^5)^{24}_{\infty}}{(q)^{25}_{\infty} }
          +  5^{12} \cdot q^4~ \frac{(q^5)^{30}_{\infty}}{(q)^{31}_{\infty} } ,
            \end{eqnarray}     
           \end{subequations}      
 where $(q)_{\infty} \equiv \prod_{m=1}^{\infty} (1-q^m)$.  The first identity, combined with eq.(26),  gives us 
           \begin{subequations}
           \begin{eqnarray}
           \frac{ p(5k+4)} {5} &=&  \left|  \begin{array} {c} ~\\~\\~\\~\\~ \end{array}  \left( \sum_{m \geq 0} (-1)^m(2m+1)J^{m(m+1)/2}
            \right)^2 \right. \nonumber  \\ 
            &&  + \left. \left( I+  \sum_{m=\pm1,\pm2, \ldots }(-1)^{m}  J^{5m(3m-1)/2} \right)^5 \times \left( \left.
             \begin{array}{cccc}
            ~ & ~  & ~&~  \\
            ~ & ~  & ~&~  \\
            ~ &  \mbox{{\rm \Huge 0}} & ~ &~ \\
            ~ & ~ & ~ &~\\
            ~ &  ~ & ~&~  \end{array} \right|
            \begin{array}{c} 
             ~1~\\
             ~0~\\
             \vdots\\
             ~0~ \\
             - 1  \end{array} \right)    \right|_{ (k+1) \times (k+1)} \nonumber \\
              &=&  \left| \begin{array}{cccccccccc}
             ~1& ~ & ~&~&~&~&~&~&~&~1\\
            -6& ~1& ~ & ~&~&~&~&~&~&~0\\
            ~9& -6& ~1& ~ & ~&~&~&~&~&~0\\
            ~10& ~9& -6& ~1 & ~ & ~&~&~&~&~0\\
            -30 & ~10 & ~9& -6&~1 & ~ &~&~& ~&~0\\
            ~0& -30 & ~10& ~9&   -6&~1& ~ &~&~& -5\\
            ~11&  ~0& -30 & ~10& ~9& -6&~1 &  ~ &~&~0\\
             ~42& ~11&  ~0& -30 & ~10& ~9& -6&~  &  ~ &~0\\
           ~ \vdots & ~&~&~&~&~&~&~~\ddots &~&  ~\vdots 
            \end{array} \right|_{(k+1) \times (k+1) },
            \end{eqnarray}
while the second one gives
          \begin{eqnarray}
          \frac{ p(7k+5)} {7} &=&   \left|   \begin{array} {c} ~\\~\\~\\~\\~ \end{array}   \left(I  +  \sum_{m=\pm1,\pm2, \ldots  } (-1)^{m} 
           J^{m(3m-1)/2}  \right)^8  \right. \nonumber \\
          && + \left\{ \left( \sum_{m \geq 0} (-1)^m(2m+1)J^{7m(m+1)/2} \right)  \left(I  + \sum_{m=\pm1,\pm2, \ldots  }
           (-1)^{m} J^{m(3m-1)/2}  \right)^4   \right.  \nonumber \\  && +  \left. \left.  7J  \left(I  + \sum_{m=\pm1,\pm2, \ldots  } (-1)^{m} 
           J^{7m(3m-1)/2}  \right)^7  \right\} 
           \times \left( \left. \begin{array}{cccc}
            ~ & ~  & ~&~  \\
            ~ & ~  & ~&~  \\
            ~ &  \mbox{{\rm \Huge 0}} & ~ &~ \\
            ~ & ~ & ~ &~\\
            ~ &  ~ & ~&~  \end{array} \right|
           \begin{array}{c} 
             ~1~\\
             ~0~\\
             \vdots\\
             ~0~ \\
            - 1  \end{array} \right)    \right|_{ (k+1) \times (k+1)} \nonumber \\
             &=&  \left| \begin{array}{cccccccccc}
             ~1& ~ & ~&~&~&~&~&~&~&~1\\
            -8& ~1& ~ & ~&~&~&~&~&~&~3\\
            ~20& -8& ~1& ~ & ~&~&~&~&~&~2\\
            ~0& ~20& -8& ~1 & ~ & ~&~&~&~&~8\\
            -70 & ~0 & ~20& -8&~1 & ~ &~&~& ~&-5\\
            ~64& -70 & ~0& ~20&   -8&~1& ~ &~&~& -4\\
            ~56&  ~64& -70 & ~0& ~20& -8&~1  &  ~ &~&-10\\
             ~ 0&~56&  ~64& -70 & ~0& ~20& -8&~ &  ~ &~5\\   
             ~\vdots & ~&~&~&~&~&~&~~\ddots &~& ~ \vdots 
            \end{array} \right|_{(k+1) \times (k+1)}
           \end{eqnarray}  
           \end{subequations} 
(where we've used the Jacobi identity $(1-z-z^2+z^5+z^7  - \cdots)^3= 1-3z+5z^3-7z^6+9z^{10} + \cdots $).
The matrices above thus consist of an LTT ``base'' part and a ``tower'' part.  The coefficients of the base matrix for 
$p(5k+4)$ correspond to sequence A000729, the coefficients in the expansion of $\left( \prod_k (1-x^k)\right)^6$, while the coefficients of powers of $J^5$ in the tower part is sequence A000728.  Likewise,  sequence A000731 gives the base 
matrix for $p(7k+5)$, while the tower part involves a combination of sequences A000727, A000730, and A010816.

From the 3rd identity, the coefficients for the base matrix for $p(25k+24)$ are given by the  expansion 
           \begin{eqnarray}
           (q)^{31}_{\infty} &=& 1 -31q+434q^2-3565q^3 
           +18445q^4 -  57505q^5 +70091q^6 + 227447 q^7 + \cdots ,
           \end{eqnarray}
(sequence A010836);  the tower part is from the expansion 
           \begin{eqnarray}          
          && 63 ~ (q)^{24}_{\infty}  (q^5)^{6}_{\infty} + 5^3 \cdot 52q  ~(q)^{18}_{\infty}  (q^5)^{12}_{\infty}  
          +5^5 \cdot 63q^2 ~(q)^{12}_{\infty}   (q^5)^{18}_{\infty} \nonumber \\ && ~~~~+ 5^{8} \cdot 6q^3 ~(q)^{6}_{\infty} 
            (q^5)^{24}_{\infty} +  5^{10} \cdot q^4~ (q^5)^{30}_{\infty}  \\
           && ~~~~~= 63 + 4988q + 95751q^2 +766014q^3 + 3323665q^4 + 8359848q^5 \nonumber \\
           && ~~~~~~+ 10896075 q^6  -6659766 q^7 +  \cdots ,\nonumber
          \end{eqnarray}
and so we have,  to this order,      
          \begin{eqnarray}
          \frac{p(25k+24)}{25} = \left| \begin{array}{cccccccccc}
          ~1& ~ & ~&~&~&~&~&~&~&~63\\
          -31& ~1& ~ & ~&~&~&~&~&~&~4988\\
         ~434& -31& ~1& ~ & ~&~&~&~&~&~95751\\
         -3565& ~434& -31& ~1 & ~ & ~&~&~&~&~766014\\
         ~18445 & -3565 & ~434& -31&~1 & ~ &~&~& ~&~3323665\\
         -57505& ~18445 & -3565& ~434&   -31&~1&  ~ &~& ~&~8359848\\
         ~70091&  -57505& ~18445 & -3565& ~434& -31&~1 &  ~ &~&~10896075\\
        ~ 227447 & ~70091&  -57505& ~18445 & -3565& ~434& -31& ~ &~& -6659766~\\    
         ~\vdots & ~&~&~&~&~& ~&~~\ddots &~&  ~\vdots 
        \end{array} \right|_{ (k+1) \times (k+1)}
        \end{eqnarray}
  We list below some sample calculations of partition functions using these determinants:
     \begin{eqnarray*}
      p(24) &=&5 \cdot  \left| \begin{array}{ccccc} ~1 &~&~&~&~1 \\ -6 & ~1& ~&~&~0\\ ~9 & -6& ~1& ~&~0\\~10&~9&-6&~1&~0\\
      -30&~10&~9& -6&~0 \end{array} \right| = 1575; \\ ~\\       
      p(40) &=&  7 \cdot   \left|  \begin{array}{cccccc} ~1 &~&~&~&~&~1 \\ -8 & ~1& ~&~&~&~3\\ ~20 & -8& ~1& ~&~&~2\\
      ~0&~20&-8&~1&~& ~8\\-70&~0&~20& -8&~1&-5\\ ~64&-70&~0&~20& -8&-4   
               \end{array}      \right|  = 37338; \\ ~\\
           p(199) &=& 25 \cdot \left| \begin{array}{cccccccc}
         ~1& ~ & ~&~&~&~&~&63\\
         -31& ~1& ~ & ~&~&~&~&4988\\
        ~434& -31& ~1& ~ & ~&~&~&95751\\
        -3565& ~434& -31& ~1 & ~ & ~&~&766014\\
        ~18445 & -3565 & ~434& -31&~1 & ~& ~&3323665\\
        -57505& ~18445 & -3565& ~434&   -31&~1&  ~&8359848\\
        ~70091&  -57505& ~18445 & -3565& ~434& -31& 1 &~10896075\\
        ~227447 & ~70091&  -57505& ~18445 & -3565& ~434& -31 &-6659766\\
        \end{array} \right| \\~\\ &=&   3646072432125.
            \end{eqnarray*}

   \section{ $p(5k+a)$ and $p(25k+a)$ determinants}     
 One can ask if it's possible to fill in some of the gaps in equations (A2a), (A2b) and (A5) and to get expressions for $p(5k+a)$,
 etc., for other values of $a$.  In the following we will consider the problem of  generalizing eqs.  (A2a) and (A5).  Ramanujan 
 \cite {Ram} derived the relation    
       \begin{eqnarray}
    \frac{ (q^5)_{\infty}}{(q^{1/5})_{\infty}} 
       &=&  \frac{(J_1^4+3qJ_2) +q^{1/5}(J_1^3+2qJ_2^2) +q^{2/5} (2J_1^2+qJ_2^3)
       +q^{3/5}(3J_1+qJ_2^4) +5q^{4/5}}{J_1^5 -11q + q^2J_2^5}  
              \end{eqnarray}  
 (his eq.(20.5)), where  he defined the  functions $J_1(q)$ and $J_2(q)$ by the equation 
              \begin{eqnarray}
      \frac{ (q^{1/5})_{\infty}} {(q^5)_{\infty}}= J_1-q^{1/5}+q^{2/5}J_2.
       \end{eqnarray}
$J_1$ and $J_2$ are series expansions in $q$ with integer coefficients and exponents.  Ramanujan then proved the
identities
     \begin{eqnarray}
    J_1^5 -11q + q^2J_2^5
         &=&  \frac {(q)^6_{\infty}  } { (q^5)^6_{\infty}}; ~~J_1J_2 =-1.
       \end{eqnarray}  
($J_1(q)$ and $J_2(q)$  are given by sequences A003823 and A007325, respectively.)  From the first identity and 
from (B1) we have  
       \begin{eqnarray}
       \sum_{n=0}^{\infty} p(n) q^{n/5}  &=&   \left[ (J_1^4+3qJ_2) +q^{1/5}(J_1^3+2qJ_2^2) +q^{2/5} (2J_1^2+qJ_2^3)\right. 
       \nonumber \\
       && ~~~~~~~~~~ \left. +q^{3/5}(3J_1+qJ_2^4) +5q^{4/5}\right] ~\frac{ (q^5)^5_{\infty}}{(q)^6_{\infty}},
        \end{eqnarray} 
  and so       
       \begin{eqnarray*}
       \sum_{k=0}^{\infty} p(5k) q^{k}  &=&  (J_1^4+3qJ_2)~\frac{ (q^5)^5_{\infty}}{(q)^6_{\infty}},\\
        \sum_{k=0}^{\infty} p(5k+1) q^{k}  &=& (J_1^3+2qJ_2^2)  ~\frac{ (q^5)^5_{\infty}}{(q)^6_{\infty}},\\
       \sum_{k=0}^{\infty} p(5k+2) q^{k}  &=&  (2J_1^2+qJ_2^3) ~\frac{ (q^5)^5_{\infty}}{(q)^6_{\infty}},\\
        \sum_{k=0}^{\infty} p(5k+3) q^{k}  &=& (3J_1+qJ_2^4)  ~\frac{ (q^5)^5_{\infty}}{(q)^6_{\infty}},     
         \end{eqnarray*} 
in addition to (A1a).   Let
          \begin{eqnarray}
         (q^{1/5})_{\infty} =  G_1 -q^{1/5}(q^5)_{\infty} + q^{2/5} G_2 ;   ~~ J_{1,2} = \frac{ G_{1,2}}{(q^5)_{\infty}}.
          \end{eqnarray}
Then, for $a=0,1,2,3,4$,                       
       \begin{eqnarray}
        \sum_{k=0}^{\infty} p(5k+a) q^{k}  &=&~\frac{ 1}{(q)^6_{\infty}} \left[ ~F_{a+1}(q^5)^{a+1}_{\infty} G_1^{4-a}
        +F_{4-a} q(q^5)^{4-a}_{\infty} G_2^{a+1} ~\right]   
         \end{eqnarray} 
 where $F_n$ is the $n$-th Fibonacci~number,  with $F_0 =0$.   We then have, making the replacement
 $q \rightarrow J$,
       \begin{eqnarray}
        p(5k+a)  &=& \left| ~(J)^6_{\infty} + \left[~F_{a+1}(J^5)^{a+1}_{\infty} G_1^{4-a}(J)
        +F_{4-a} J(J^5)^{4-a}_{\infty} G_2^{a+1}(J)~ \right]  
         \times \left( \left.
             \begin{array}{cccc}
            ~ & ~  & ~&~  \\
            ~ & ~  & ~&~  \\
            ~ &  \mbox{{\rm \Huge 0}} & ~ &~ \\
            ~ & ~ & ~ &~\\
            ~ &  ~ & ~&~  \end{array} \right|
            \begin{array}{c} 
             ~1~\\
             ~0~\\
             \vdots\\
             ~0~ \\
             - 1  \end{array} \right) \right|_{ (k+1)\times(k+1)} \nonumber \\
             &=&  \left| \begin{array}{cccccccccc}
             ~1& ~ & ~&~&~&~&~&~&~&X_0\\
            -6& ~1& ~ & ~&~&~&~&~&~&X_1\\
            ~9& -6& ~1& ~ & ~&~&~&~&~&X_2\\
            ~10& ~9& -6& ~1 & ~ & ~&~&~&~&X_3\\
            -30 & ~10 & ~9& -6&~1 & ~ &~&~& ~&X_4\\
            ~0& -30 & ~10& ~9&   -6&~1& ~ &~&~& X_5\\
            ~11&  ~0& -30 & ~10& ~9& -6&~1 &  ~ &~&X_6\\
             ~42& ~11&  ~0& -30 & ~10& ~9& -6&~ &  ~ &X_7\\
           ~ \vdots & ~&~&~&~&~&~&~~\ddots &~&  ~\vdots 
            \end{array} \right|_{(k+1) \times (k+1)},
         \end{eqnarray}   
where the elements $X_n(=X_n^{(a)})$ in the tower matrix are determined for each value of $a$ by the expansion 
          \begin{eqnarray}
          F_{a+1}(J^5)^{a+1}_{\infty} G_1^{4-a}(J) +F_{4-a} J(J^5)^{4-a}_{\infty} G_2^{a+1}(J)= X_0I +X_1J +X_2J^2 + \cdots
          \end{eqnarray}

The $G$'s can be expressed in terms of the Ramanujan theta function:
   \begin{eqnarray}
   G_1(q) = \frac{f(-q^2,-q^3)^2}{ f(-q,-q^2)} ;~~
   G_2(q) = -\frac{f(-q,-q^4)^2}{ f(-q,-q^2)} ;~~
        f(a,b) = \sum_{n= -\infty}^{\infty} a^{n(n+1)/2}b^{n(n-1)/2}.
      \end{eqnarray}
      However, it is more convenient to write them as the series expansions
      \begin{eqnarray*}
      G_1(q) &=& -1 + \sum_{k=0}^{\infty} q^{k(30k-1)} \left[ 1 + q^{2k} +q^{12k+1} -q^{20k+3} -q^{30k+7}
          -q^{32k+8} -q^{42k+14} + q^{50k+20}  \right]\\
          &=&  1 +q-q^3   -q^7 -q^8-q^{14} + q^{20} + q^{29} + q^{31} + q^{42} -  q^{52} + \cdots, \\
           G_2(q) &=&  \sum_{k=0}^{\infty} q^{k(30k+7)} \left[ -1 + q^{6k+1} -q^{10k+2} +q^{16k+4} +q^{30k+11}
          -q^{36k+15} +q^{40k+18} - q^{46k+23}  \right] \\
          &=&  -1 +q-q^2  +q^4+q^{11} -q^{15} + q^{18} - q^{23} - q^{37}+q^{44} -q^{49} +q^{57} + \cdots ,
      \end{eqnarray*}
which follow directly from (B5).  ($G_1(q) $ corresponds to sequence A113681, and $-G_2(q)$ to sequence
A116915.)  In column-vector form, the first terms in the expansions from eq. (B8) are:
\begin{eqnarray*}
{\bf X}^{(0)}  = \left( \begin{array}{c} ~~1 \\ ~~1 \\ ~~9 \\ -3 \\ -11 \\ -10 \\ ~~10 \\ -10 \\   \vdots~ \end{array} \right);
~~{\bf X}^{(1)}  = \left( \begin{array}{c} ~~1 \\ ~~5  \\ -1 \\ ~~4 \\ -10 \\ -7 \\-5 \\ ~~2 \\ \vdots ~\end{array} \right);
~~{\bf X}^{(2)} = \left( \begin{array}{c} ~~2 \\ ~~3  \\ ~~5 \\ -10 \\ ~~3 \\ -9 \\-11 \\ -8  \\ \vdots ~\end{array} \right);
~~{\bf X}^{(3)}  = \left( \begin{array}{c} ~~3 \\ ~~4  \\ -4 \\ ~~7 \\ -16 \\ ~~3 \\ -17 \\ -13  \\ \vdots ~\end{array} \right);
\end{eqnarray*}
with the coefficients for ${\bf X}^{(4)}$ being given in (A2a).

We can generalize eq.(A5) to an expression for $p(25k +a)$ for $a=4,9,14$ and 19 following Ramanujan's derivation of the identity (A1c)\cite{Ram}.  We make the replacement $q \rightarrow q^{1/5}$ in eq. (A1a) and get
           \begin{eqnarray}
          \sum_{k=0}^{\infty} p(5k+4)q^{k/5} = 5~ \frac{(q)^5_{\infty}}{(q^{1/5})^6_{\infty} } 
             =5 ~ \frac{(q)^5_{\infty}}{ (q^5)_{\infty}^6 } ~\frac{1} { (J_1-q^{1/5}+q^{2/5}J_2)^6}  .
           \end{eqnarray}
To simplify the notation, we define $ x \equiv q/J_1^5=-qJ_2^5 $.    Then
          \begin{eqnarray}
            J_1-q^{1/5}+q^{2/5}J_2= J_1 (1-x^{1/5}-x^{2/5}) = J_1 (1+x^{1/5}/\phi )(1-\phi~x^{1/5})
            \end{eqnarray}
where $\phi = ( \sqrt{5} +1)/2$, the golden ratio.    We have that                 
           \begin{eqnarray}
           \frac{1}{1+ ax^{1/5} } = \frac{1- ax^{1/5} + a^2x^{2/5} - a^3x^{3/5} +a^4x^{4/5}} {1+ a^5x} .
           \end{eqnarray}
Then           
             \begin{eqnarray}
          \frac{1} { 1-x^{1/5}-x^{2/5}}&=& \frac{1} { (1+x^{1/5}/\phi)(1-\phi x^{1/5})} \nonumber \\ &=&   
           \frac{1-x^{1/5}/\phi + x^{2/5}/\phi^2- x^{3/5}/\phi^3 +x^{4/5}\/\phi^4} {1+x/\phi^5}
             ~ \frac{1+\phi x^{1/5} + \phi^2x^{2/5}+\phi^3x^{3/5} +\phi^4x^{4/5}} {1-\phi^5x} \nonumber \\
             &=& \frac{    1-3x + (1+2x) x^{1/5}+  (2-x)x^{2/5}+ (3+x) x^{3/5}+5x^{4/5}    } {  1-11x-x^2  }  
            \end{eqnarray}
The denominator in (B13) is            
             \begin{eqnarray}
             1-11x-x^2  = \frac{J_1^5 -11q +q^2J_2^5 } {J_1^5} =  \frac{1}{J_1^5} \frac {(q)^{6}_{\infty}  } { (q^5)^6_{\infty}} 
             \end{eqnarray}
The numerator in (B13) is to be raised to the sixth power  in eq.(B10).  We define the functions $H_n(x)$ to be series expansions in $x$ with integer coefficients and exponents such that
         \begin{eqnarray}
        \left[  1-3x + (1+2x) x^{1/5}+  (2-x)x^{2/5}+ (3+x) x^{3/5}+5x^{4/5} \right]^6 &=& H_1(x) +H_2(x)x^{1/5} 
        + H_3(x)x^{2/5} \nonumber \\
         && ~~~ + H_4(x) x^{3/5} + H_5(x) x^{4/5}
        \end{eqnarray}    
Expanding the left side of this equation and collecting terms, we get
 \begin{eqnarray*}
          H_1(x) &=& -98 x^9 +9939 x^8 -107712 x^7 + 167031 x^6 - 27918 x^5 + 127011 x^4 + 160552 x^3 +
               32784 x^2 + 858 x+1\\
          H_2(x) &=&  ~~27 x^9 -4806 x^8 + 78758 x^7 - 171984 x^6 +98667 x^5 +78986 x^4 + 176592 x^3+52644 x^2 
          + 2138 x +6\\ 
          H_3(x) &=& \; -6 x^9 +2138 x^8 -52644 x^7 + 176592 x^6 - 78986 x^5 + 98667 x^4 + 171984 x^3 +
               78758 x^2 + 4806 x+27 \\  
         H_4(x) &=& ~~~~ x^9 - 858 x^8 + 32784 x^7  - 160552 x^6 + 127011 x^5 + 27918 x^4 + 167031 x^3 +
               107712  x^2 + 9939 x+98\\                                                                   
          H_5(x) &=&  ~~~~~~~~~~~~\, 315 x^8 -18640 x^7 + 139305 x^6 -127020 x^5 + 106425 x^4 + 127020 x^3 +
               139305 x^2 + 18640 x+315
          \end{eqnarray*}    
Then, changing back to the variable $q$, eq.(B10) is    
      \begin{eqnarray}
       \sum_{k=0}^{\infty} p(5k+4)q^{k/5} &=& 5~ \frac{(q^5)^{30}_{\infty}}{ (q)^{31}_{\infty}  }~\left( J_1^{24} H_1(q) 
       + J_1^{23}H_2(q)q^{1/5}+ J_1^{22} H_3(q)q^{2/5}  +J_1^{21}H_4(q) q^{3/5}+ J_1^{20} H_5(q) q^{4/5}\right)        
     \end{eqnarray}
and we have   
      \begin{eqnarray*}
      \sum_{k=0}^{\infty} p(25k+4)q^{k} &=&  5~ \frac{(q^5)^{30}_{\infty}}{ (q)^{31}_{\infty}  } 
       \left[ ~ J_1^{24} + 858~ qJ_1^{19}+  32784~ q^2J_1^{14}   +160552~ q^3J_1^{9} + 127011~ q^4J_1^4  \right. \\ 
       && \left.~~~~~~~~~~~~  + 27918~ q^5J_2+ 167031 ~q^6J_2^{6} +107712~q^7J_2^{11}+9939~ q^8J_2^{16}
       +98~ q^9J_2^{21} \right] \\~ \\
       \sum_{k=0}^{\infty} p(25k+9)q^{k} &=&  5~ \frac{(q^5)^{30}_{\infty}}{ (q)^{31}_{\infty}  }  
        \left[ ~ 6J_1^{23} + 2138~ qJ_1^{18}+ 52644~ q^2J_1^{13}   +176592~ q^3J_1^{8} + 78986~ q^4J_1^3 \right. \\ 
        && \left.~~~~~~~~~~~~ + 98667~ q^5J_2^2+ 171984 ~q^6J_2^{7} +78758~q^7J_2^{12}
        +4806~ q^8J_2^{17}+ 27~q^9J_2^{22}   \right] \\ ~\\      
        \sum_{k=0}^{\infty} p(25k+14)q^{k}  &=& 5~ \frac{(q^5)^{30}_{\infty}}{ (q)^{31}_{\infty}  } 
        \left[ ~ 27J_1^{22} + 4806~ qJ_1^{17}+ 78758~ q^2J_1^{12}   +171984~ q^3J_1^{7} + 98667~ q^4J_1^2 \right. \\ 
         && \left.~~~~~~~~~~~~ + 78986~ q^5J_2^3+ 176592 ~q^6J_2^{8} +52644~q^7J_2^{13}
          +2138~ q^8J_2^{18}+6~ q^9J_2^{23} \right] \\ ~\\   
      \sum_{k=0}^{\infty} p(25k+19)q^{k} &=& 5~ \frac{(q^5)^{30}_{\infty}}{ (q)^{31}_{\infty}  } 
      \left[~ 98J_1^{21}+ 9939 ~qJ_1^{16}  + 107712 ~ q^2J_1^{11} + 167031~ q^3J_1^{6}  + 27918 ~q^4 J_1   \right. \\
          &&\left.  ~~~~~~~~~~~~   + 127011~ q^5J_2^{4} +160552~ q^6J_2^{9}
             + 32784 ~q^7J_2^{14}+858 ~q^8J_2^{19} + q^9J_2^{24} \right]     \\ ~ \\               
          \sum_{k=0}^{\infty} p(25k+24)q^{k} &=& 5^2 \frac{(q^5)^{30}_{\infty}}{ (q)^{31}_{\infty}  }
              \left[ ~ 63 J_1^{20} +3728 ~qJ_1^{15}+27861~ q^2J_1^{10}+25404~ q^3J_1^{5}  + 21285~ q^4\right. \\
               && ~~~~~~~~~~~~~ \left. +25404 ~q^5J_2^{5}+ 27861~ q^6J_2^{10} 
               +3728~ q^7J_2^{15}  + 63 ~q^8 J_2^{20}\right]                     
      \end{eqnarray*}
 The expression on the right in the last equation above reduces to Ramanujan's result in (A1c) upon the substitutions
       \begin{eqnarray*}
        J_1^5 +q^2J_2^5 &=& X +11~q, \\
        J_1^{10} +q^4J_2^{10}  &=&  X^2 +22~q X +123~q^2,\\
        J_1^{15} +q^6J_2^{15}  &=& X^3 +33~qX^2 +366~q^2 X +1364~q^3,\\
        J_1^{20} +q^6J_2^{20} &=& X^4 +44~q X^3 + 730~q^2X^2 +5412~q^3 X +15127~q^4,  
        \end{eqnarray*}
where
        \begin{eqnarray*}
        X &\equiv&  \frac {(q)^{6}_{\infty}  } { (q^5)^6_{\infty}} .
        \end{eqnarray*}  
As before, the $Z_n$ coefficients in the $p(25k+a)$ determinant     
     \begin{eqnarray}
     p(25k+a) = 5 \cdot \left| \begin{array}{cccccccccc}
          ~1& ~ & ~&~&~&~&~&~&~&Z_0\\
          -31& ~1& ~ & ~&~&~&~&~&~& Z_1\\
         ~434& -31& ~1& ~ & ~&~&~&~&~& Z_2\\
         -3565& ~434& -31& ~1 & ~ & ~&~&~&~& Z_3\\
         ~18445 & -3565 & ~434& -31&~1 & ~ &~&~& ~& Z_4\\
         -57505& ~18445 & -3565& ~434&   -31&~1&  ~ &~& ~& Z_5\\
         ~70091&  -57505& ~18445 & -3565& ~434& -31&~1 &  ~ &~& Z_6\\
        ~ 227447 & ~70091&  -57505& ~18445 & -3565& ~434& -31& ~ &~& Z_7 \\    
         ~\vdots & ~&~&~&~&~& ~&~~\ddots &~&  ~\vdots 
        \end{array} \right|_{ (k+1) \times (k+1)}
        \end{eqnarray}
are obtained by an expansion in powers of $q$ of the numerators on the RHS 's of these generating-function equations.
We have, for $a=4,9,14, 19$,   
   \begin{eqnarray*}
{\bf Z}^{(4)}  = \left( \begin{array}{c} 1 \\ 882 \\ 49362 \\ 768246 \\  5380497 \\ 20802996 \\ 47413915 \\ 46923084 \\  
 \vdots~ \end{array} \right);
~~{\bf Z}^{(9)}  = \left( \begin{array}{c} 6 \\ 2276  \\ 92646 \\ 1198566 \\ 7354172 \\ 25710039 \\ 51224670 \\ 39450895 \\ 
\vdots ~\end{array} \right);
~~{\bf Z}^{(14)} = \left( \begin{array}{c}  27 \\ 5400 \\ 166697 \\ 1811682 \\ 9871992 \\ 30828786 \\ 55015749 \\
        20079168 \\ \vdots ~\end{array} \right);
~~{\bf Z}^{(19)}  = \left( \begin{array}{c} 98 \\ 11997  \\ 287316 \\ 2672825 \\12906450 \\   36553962 \\  54917174
     \\ 2443563\\  \vdots ~\end{array} \right).
\end{eqnarray*}

  \section{$\sum p(n)x^{n}$ determinants }
             
Expression (11) for $p(n)$ can be used to express finite sums of the form $\sum p(n)x^{n}$ as determinants.  We have
     \begin{eqnarray}
     \sum_{n=0}^k p(n) x^{n} = x^{k} \sum_{n=0}^k x^{n-k}p(n) 
       = x^{k} \left| \begin{array} {cccccccc}
       ~~1 & ~ & ~ & ~ &~&~ &~&~~1 \\
       -1 & ~~1  &~&~&~&~&~& ~~1/x\\
        -1 & -1  & ~~1 & ~ &~&~&~& ~~1/x^2\\
         ~0 & -1  &-1&~~ \ddots &~&~&~& ~~1/x^3\\       
        ~\vdots  & ~  &~ &~~\ddots & ~&~&~& ~~\vdots \\
         d_{k-1}& d_{k-2} &~&~&~&~& ~~1 &~~1/x^{k-1} \\
        d_k & d_{k-1} & \cdots &~&~&~& -1 & ~~1/x^k
                  \end{array} \right|_{(k+1)\times (k+1)} 
          \end{eqnarray}
The tower part of the matrix can be expressed as 
      \begin{eqnarray}
       (I-  J/x)^{-1} \left( \left. \begin{array}{cccc}
            ~ & ~  & ~&~  \\
            ~ & ~  & ~&~  \\
            ~ & ~  & ~&~  \\
           ~ &    \mbox{{\rm \Huge 0}} &~ &~ \\
            ~ & ~ & ~ &~\\
            ~ & ~  & ~&~  \\
            ~ &  ~ & ~&~  \end{array} \right|
             \begin{array}{c} 
             ~1~\\
             ~0 ~\\
             ~0~\\
              ~\vdots ~\\
              ~0~\\
                0 \\
                -1 \end{array} \right)_{(k+1)\times (k+1)} . \nonumber
              \end{eqnarray}   
We now multiply the determinant in eq. (C1) by the determinant of $(I-J/x)$, (which is equal to 1), and get              
              \begin{eqnarray}
      \sum_{n=0}^k p(n) x^{n}&=& x^k   \left|  \begin{array} {cccccccc}
          ~~1 & ~ & ~ & ~ &~&~ &~&~~1 \\
       -1-1/x  & ~~1  &~&~&~&~&~& ~~0\\
        -1+1/x  & -1-1/x  & ~~1 & ~ &~&~&~& ~~0\\
         ~1/x & -1+1/x  &-1-1/x &~~ \ddots &~&~&~& ~~0\\       
        ~\vdots  & ~  &~ &~~\ddots & ~&~&~& ~~\vdots \\
         d_{k-1}-d_{k-2}/x & d_{k-2}-d_{k-3}/x  &~&~&~&~& ~~1 &~~0 \\
         d_{k}-d_{k-1}/x & d_{k-1} -d_{k-2}/x & \cdots &~&~&~& -1-1/x & ~~0
                  \end{array}       \right|_{(k+1)\times (k+1)}   \nonumber \\ \nonumber ~\\
           &=&  \left|  \begin{array} {cccccccc}
          x+1 & -x & ~ & ~ &~&~ &~&~ \\
        x-1  & x+1  &-x &~&~&~&~& ~\\
         -1  & x-1  & x+1 &  ~~~~\ddots &~&~&~& ~\\
         0 & -1  & x-1 &~~~~ \ddots &~&~&~& ~\\       
        \vdots  & ~  &~ &~~~~\ddots & ~&~&~& ~ \\
         d_{k-2}-xd_{k-1} & d_{k-3}-xd_{k-2}  &~&\cdots &~&~& x+1 &-x \\
         d_{k-1} -xd_{k} & d_{k-2} -xd_{k-1} & ~& \cdots &~&~& x-1 & x+1
                  \end{array}       \right|_{k\times k}  ,       
          \end{eqnarray}  
where in the last line we've expanded the $(k+1)$-dimensional determinant by minors along the final column and then taken the factor  $(-1)^k x^k$ inside the resulting $k$-dimensional determinant, multiplying each of the columns by $-x$.  The final result can be written in a more compact notation as               
     \begin{eqnarray}
           \sum_{n=0}^{k} p(n) x^n = \det \left[  -xJ^T +I +\sum_{m > 0}^{q_m<k+1}  (-1)^m \left[ J^{m(3m-1)/2} +J^{m(3m+1)/2} 
                  -x J^{(m-1)(3m+2)/2} -x J^{(m+1)(3m-2)/2} \right]     \right]   _{k \times k} . 
     \end{eqnarray}

 ~\\
2000 Mathematics Subject Classification:  Primary:  05A17; Secondary:  11Y35, 11B68.\\
Keywords:  partition function, Bernoulli numbers,  Euler numbers, Sterling numbers. \\
(Concerned with sequences A000009, A000041, A000364, A000728, A000729, A000730, A000731,
 A001318, A010815, A010816, A010836, A095699, A113681, and A116915)

\end{document}